# Generalized linearization of nonlinear algebraic equations: an innovative approach

W. Chen

Permanent mail address: P. O. Box 2-19-201, Jiangshu University of Science & Technology, Zhenjiang City, Jiangsu Province 212013, P. R. China

Present mail address (as a JSPS Postdoctoral Research Fellow): Apt.4, West 1st floor, Himawari-so, 316-2, Wakasato-kitaichi, Nagano-city, Nagano-ken, 380-0926, JAPAN

E-mail: chenw@homer.shinshu-u.ac.jp

Permanent email: chenwwhy@hotmail.com

## Introduction

Based on the matrix expression of general nonlinear numerical analogues [1, 2], this paper proposes a novel philosophy of nonlinear computation and analysis. The nonlinear problems are considered an ill-posed linear system. In this way, all nonlinear algebraic terms are instead expressed as Linearly independent variables. Therefore, a n-dimension nonlinear system can be expanded as a linear system of n(n+1)/2 dimension space. This introduces the possibility to applying generalized inverse of matrix to computation of nonlinear systems. Also, singular value decomposition (SVD) can be directly employed in nonlinear analysis by using such a methodology.

## Matrix expression of general nonlinear numerical formulation

In this section, by using the Hadamard product, we derive a generalized matrix formulation of nonlinear numerical analogues [1, 2].

**Definition 2.1** Let matrices $A=[a_{ij}]$ and $B=[b_{ij}] \in C^{N \times M}$, the Hadamard product of matrices is defined as $A \circ B = [a_{ij} b_{ij}] \in C^{N \times M}$. where $C^{N \times M}$ denotes the set of $N \times M$ real matrices.

**Theorem 2.1**: letting A, B and $Q \in C^{N \times M}$, then

1> $A \circ B = B \circ A$ (1a)

2> $k(A \circ B) = (kA) \circ B$, where k is a scalar. (1b)

3> $(A+B) \circ Q = A \circ Q + B \circ Q$ (1c)

4> $A \circ B = E_N^T (A \otimes B) E_M$, where matrix $E_N$ (or $E_M$) is defined as $E_N = [e_1 \otimes e_1 \vdots \cdots \vdots e_N \otimes e_N]$, $e_i = [0 \cdots 0 \underset{i}{1} 0 \cdots 0]$, i=1, ..., N, $E_N^T$ is the transpose matrix of $E_N$. $\otimes$ denotes the Kronecker product of matrices. (1d)

5> If A and B are non-negative, then $\lambda_{min}(A) \min\{b_{ii}\} \leq \lambda_j (A \circ B) \leq \lambda_{max}(A) \max\{b_{ii}\}$, where $\lambda$ is the eigenvalue. (1e)

6> $(\det A)(\det B) \leq \det(A \circ B)$, where det( ) denotes the determinant. (1f)

For more details about the Hadamard product see [3, 4].

It is well known that the majority of popular numerical methods such as the finite element, boundary element, finite difference, Galerkin, least square, collocation and spectral methods have their root on the method of weighted residuals (MWR) [5, 6]. Therefore, it will be generally significant to apply the Hadamard product to the nonlinear computation of the MWR. In the MWR, the desired function u in the differential governing equation

$$\psi\{u\} - f = 0, \quad \text{in} \quad \Omega \tag{2}$$

is replaced by a finite series approximation $\hat{u}$,

$$u = \hat{u} = \sum_{k=1}^{N} c_k \phi_k, \tag{3}$$

where $\psi\{\ \}$ is a differential operator. $\phi_k$ can be defined as the assumed functions and $c_k$'s are the unknown parameters. The approximate function $\hat{u}$ is completely specified in terms of unknown parameters $c_k$. Boundary conditions of governing equation can be specified

$$u = \bar{u}, \quad \text{on} \quad \Gamma_1 \tag{4a}$$

$$q = \partial u/\partial n = \bar{q}, \quad \text{on} \quad \Gamma_2, \tag{4b}$$

where n is the outward normal to the boundary, $\Gamma = \Gamma_1 + \Gamma_2$, and the upper bars indicate known boundary values.

Substituting this approximation $\hat{u}$ into the governing equation (2), it is in general unlikely that the equation will be exactly satisfied, namely, result in a residual R

$$\psi\{\hat{u}\} - f = R \tag{5}$$

The method of weighted residuals seeks to determine the N unknowns $c_k$ in such a way that the error R is minimized over the entire solution domain. This is accomplished by requiring that weighted average of the error vanishes over the solution domain. Choosing the weighting function $W_j$ and setting the integral of R to zero:

$$\int_D [\psi\{\hat{u}\} - f] W_j dD = \int_D R W_j dD = 0, \quad j=1,2,\ldots,N. \qquad (6)$$

Equation (6) can be used to obtain the N unknown coefficients. This equation also generally describes the method of weighted residuals. In order to expose our idea clearly, Let us consider the quadratic nonlinear operator of the form:

$$p(u)r(u) + L(u) = f, \qquad (7)$$

where p(u), r(u) and L(u) are linear differential operators, f is the constant. The general scheme of weighted residuals approximates Eq. (7) is given by

$$\int_\Omega [p(\hat{u})r(\hat{u}) + L(\hat{u}) - f]\phi_j d\Omega = \int_{\Gamma_2} (q - \bar{q})\phi_j d\Gamma - \int_{\Gamma_1} (u - \bar{u})\frac{\partial \phi_j}{\partial n} d\Gamma, \quad j=1,2,\ldots,N. \qquad (8)$$

Substitution of Eq. (3) into Eqs. (8) and applying equation (1d) in the theorem 2.1 result in

$$D_{n\times n} x + G_{n\times n^2} (x \otimes x) = b, \qquad (9)$$

where $\otimes$ denotes the Kroneker product of matrices; $x = \{c_k\}$, $c_k$'s are the undetermined parameters in Eq. (3).

$$D = \int_\Omega L\{\phi\} W_j d\Omega, \qquad (10a)$$

$$G_{n\times n^2} = \int_\Omega [p(\phi) \otimes r(\phi)] W_j d\Omega \in C^{n \times n^2}. \qquad (10b)$$

For the cubic nonlinear differential equations, we can obtain similar general matrix formulation by using the same approach:

$$D_{n\times n} x + R_{n\times n^3} (x \otimes x \otimes x) = b, \qquad (11)$$

where D and R are constant coefficient matrices. To simplify notation, formulations with form of Eqs. (9) and (11) are denoted as formulation-K, where K is chosen since the Kronecker product is used in expressing nonlinear numerical discretization term.

As was mentioned earlier, most of popular numerical techniques can be derived from the method of weighted residual. The only difference among these numerical methods lies in the use of different weighting and basis functions in the MWR. From the foregoing deduction, it is noted that Eq. (9) can be obtained no matter what weighting and basis functions we use in the method of weighted residuals. Therefore, it is straightforward that we can obtain the formulation-K for the nonlinear computations of all these methods.

**Generalized linearization of nonlinear variables**

It is a usual practice to use the Newton-Raphson method to solve nonlinear algebraic equations. However, if we consider the nonlinear system as an ill-posed linear problem, the nonlinear terms of nonlinear equations can be instead regarded as independent system variables. In this way, we can state a n-dimension nonlinear system as a n(n+1)/2 dimension linear one.

Eq. (9) expresses general numerical analogue of quadratic nonlinear problem in matrix form. Without the loss of generality, let

$$y_i = x_i, \qquad i=1,2,\ldots, n. \tag{12a}$$

$$y_i = x_i x_j, \qquad i=n+1, n+2, \ldots, n(n-1)/2; \quad j=1,2, \ldots, n. \tag{12b}$$

Eq. (9) can be restated as

$$P_{n \times (n+1)n/2}\, y = b, \tag{13}$$

where P can be easily obtained from the known coefficient matrices of D and G in equation (9). It is noted that Eq. (13) is an ill-posed linear equation system. Therefore, there exists more than one set of solutions. This is in agreement with the real nonlinear solutions.

For linear system equations like Eq. (13), we can not use the normal numerical solution technique. There are some approaches available for analysis and computation of such equations. Singular value decomposition (SVD) can be employed to analyze the structure features of equation (13). This may provide some useful information of nonlinear system. For example, some intrinsic structures of bifurcation, chaos, and shock wave may be illustrated via the SVD analysis of coefficient matrix P in equation (13).

Another useful approach is generalized inverse such as G-inverse and Moore-Penrose inverse [4, 7]. For instance, by using Moore-Penrose inverse concept, we can get only one definite generalized inverse of matrix P.

Finally, we provide a least square solution of equation (13). Let the transpose of matrix P premultiple equation (13), we have

$$P^T P y = P^T b. \tag{14}$$

The above equation (14) can be easily solved to obtain one set of solution of y. Then, substituting the solution into equations (12a, b) produce numerous sets of solutions of original variable x. The choice of the desired results from them should be mostly dependent on the physical background of the practical problem.

This study is still in a very early stage. Presentation of the above idea herein is to spark more beneficial discussions on the possible usefulness of this strategy.